\documentclass[12 pt]{amsart}
\usepackage[english]{babel}
\usepackage{amsfonts}
\usepackage[latin1]{inputenc}

\title[]
{Thermodynamic Dictionary for \\ Partitioned Dynamical Systems}

\author[Eduardo Santana]{Eduardo Santana}
\address{Universidade Federal de Alagoas, 
	Av. Beira Rio, s/n, 57200-000 Penedo, Brazil.}
\email{jemsmath@gmail.com}

\newtheorem{theorem}{Theorem}

\newtheorem{lemma}[theorem]{Lemma}
\newtheorem{proposition}[theorem]{Proposition}

\theoremstyle{definition}

\newtheorem{remark}[theorem]{Remark}

\newtheorem{theorema}{Theorem}

\date{\today}
\begin{document}
	
	\begin{abstract}
	In this paper we  establish a dictionary between the behavior of a general map $f: X \to X$ which is defined by using a decomposition into disjoint subsets $X = A \cup B$ and the thermodynamic formalism.
	
	The approach is set theoretical, topological and ergodic. We define a key topology and consider its Borel $\sigma$-algebra in order to study the thermodynamic formalism for the general setting, proving that recurrence implies periodicity, the topological entropy is zero, we can decompose the invariant measures into a general sum of ergodic ones and we establish a dictionary between the cycles and the equilibrium states of continuous potentials.
	
	All this is done for the general setting, with applications for the Syracuse maps and the Collatz map in particular.

	Finally, we give several applications beyond the arithmetic realm: algebraic maps, counting closed geodesics by Margulis' formula, infinite chessboard and stochastic lattices.
	\end{abstract}
	
	\maketitle
	\section{Introduction}
	Let us denote, respectively, $\mathbb{N}_{0} = \{0, 1, 2, 3, \cdots\}$ and $\mathbb{N} = \{1, 2, 3, \cdots\}$ the set of the natural numbers, including the element $0$ or not. The Collatz map $C: \mathbb{N} \to \mathbb{N}$, famously known due to the Collatz Conjecture is defined as
	\[
	C(x) =
	\left\{
	\begin{array}{ccc}
		3x + 1 \,\, &\text{if}& \,\, x \,\, \text{is odd} \\
		\frac{x}{2} \,\, &\text{if}& \,\, x \,\, \text{is even} \\
	\end{array}
	\right.
	\]
	For $m, n \in \mathbb{N}$ odd numbers, we can define the Syracuse family of maps as
	\[
	S_{m,n}(x) =
	\left\{
	\begin{array}{ccc}
		mx + n \,\, &\text{if}& \,\, x \,\, \text{is odd} \\
		\frac{x}{2} \,\, &\text{if}& \,\, x \,\, \text{is even} \\
	\end{array}
	\right.
	\]
	We can see that $C = S_{3,1}$. Also there is a decomposition $\mathbb{N} = O \cup E$ where $O$ is the set of the odd numbers and $E$ is the set of even numbers. Also, we can consider two functions $a : O \to E$ and $b: E \to \mathbb{N}$ where $a(x) = mx + n, b(x) = \frac{x}{2}$. With this remark, we have
	\[
	S_{m,n}(x) =
	\left\{
	\begin{array}{ccc}
		a(x) \,\, &\text{if}& \,\, x \,\, x \in O \\
		b(x)\,\, &\text{if}& \,\, x \in E \\
	\end{array}
	\right.
	\]
	Inpired by this remark, we can generalize it as follows. Let $X$ be an arbitrary set and $A, B \subset X$ such that $X = A \cup B, A \cap B = \emptyset$. Consider functions $a : A \to B$ and $b : B \to X$ such that $b(B) \not \subset A, b(B) \not \subset B$. Define a map $f : X \to X$ as follows
	\[
	f(x) =
	\left\{
	\begin{array}{ccc}
		a(x) \,\, &\text{if}& \,\, x \in A \\
		b(x) \,\, &\text{if}& \,\, x \in B \\
	\end{array}
	\right.
	\] 
	We will always refer to this type of map as in \textbf{general form}. We define a \textbf{key topology} $\mathcal{T}$ generated by the collection
	\[
	\{ \{x, y\}  \mid x \in X, y \in f^{-1}(x)\}.
	\]
	Its Borel $\sigma$-algebra will be denoted by $\Sigma$.
	
	\section{Main results}
	For the general setting with respect to the function $f: X \to X$ as defined above in the Introduction, we can prove the following results
	
	\begin{theorema}\label{A}
		We have the following topological and ergodic facts with respect to $f$
		\begin{itemize}
			\item Under the key topology $\mathcal{T}$ and its Borel $\sigma$-algebra every recurrent point is periodic and the ergodic probabilities are the ones supported on the periodic orbits.
			
			\item For every $f$-invariant probability $\mu$, we have that $\mu$-almost every point is periodic.
					
			\item Every forward orbit is an open subset. 
			
			
			\item Every $f$-invariant probability  has zero entropy.
		\end{itemize}
	\end{theorema}
	
	We can establish the following dictionary between the cycles and the equilibrium states.
	
	\begin{theorema}\label{B}
		For the general map $f$ and we have the following facts related to equilibrium states.
		\begin{itemize}		
			\item The existence of an equilibrium state for every continuous potential $\phi : \mathbb{N} \to \mathbb{N}_{0}$ with finite pressure guarantees either finiteness of cycles or infinitely many cycles sharing the same period.
			
			\item The existence of a continuous potential $\phi : \mathbb{N} \to \mathbb{N}_{0}$ with finite pressure and no equilibrium state is equivalent to the unboundedness of the periods of the infinitely many cycles.
			
			\item Uniqueness of periodic orbits is guaranteed by the uniqueness of equilibrium state for every bounded and continuous potential $\phi : \mathbb{N} \to \mathbb{N}_{0}$.
		\end{itemize}
	\end{theorema}
	
	Consider the fundamental axiom on the itineraries for maps taken at random.
	
	\textbf{Fundamental Axiom on the Maps (FAM):} Let us a $f: X \to X$ of the general form. Since the set $X$ admits the following decomposition $X = A \cup B$, we have the itinerary up time $p$ of a point $x$ as $\sigma(x) \in \{A, B\}^{p}$. We assume the the following equation has at most finitely many solutions $M_{\sigma}$:
	\[
	f_{\sigma(x)} (x) = f_{\sigma_{p}(x)} \circ \cdots f_{\sigma_{1}(x)} (x) = x,
	\]
	where every map $f_{\sigma_{i}(x)} \in \{a(x), b(x)\}$. Under this hypothesis, we can establish an upper bound for the cardinality of cycles.	

	\begin{theorema}\label{C}
		For map $f_{\sigma(x)}$ satisfying (FAM) we have the following facts.
		\begin{itemize}		
			\item For each natural number $p \in \mathbb{N}$ we have at most $M_{\sigma} 2^{p}/p$ cycles with period $p$.
			
			\item There are at most countably many cycles.
			
			
						
		\end{itemize}
	\end{theorema}
	
	In the seminal work by Terence Tao \cite{Tao2022} it is established that almost every natural number reaches a low value, almost proving the conjecture. However, this work does not deal with the cardinality of cycles and the general Syracuse maps either. Also, the divergence of orbits is not discussed.
	
	The main novelty of our work is obtaining relation between the cardinality of cycles of every Syracuse map and the existence of equilibrium states (in particular for the Collatz map). It opens the way to study the Collatz Conjecture under the language of the thermodynamic formalism and questions on the Syracuse maps as well.
	
	We will prove the Theorems \ref{A}, \ref{B}, \ref{C} for the general map $f$ and so to the family of Syracuse maps $f$, having the Collatz and the Baker maps as examples. The Baker map sends every odd number $x$ into $3x - 1$, as a variation of the Collatz map.
	\section{Proof of Theorem \ref{A}}
	We divide the proof of our first main theorem into some lemmas and remarks. We proceed with this now. 
	
	\subsection{Trapping the orbits towards periodicity}
	The following lemma proves that recurrence implies periodicity.
	\begin{lemma}\label{periodic}
		Consider the topology $\mathcal{T}$ defined in the Introduction and its Borel $\sigma$-algebra $\Sigma$ with respect to the general map $f: X \to X$ possessing some $f$-invariant Borel probability $\mu$.Then $\mu$-almost every point is periodic.
	\end{lemma}
	\begin{proof}		
		We will prove that every recurrent point is indeed periodic. By Poincaré Recurrence Theorem we have that $\mu$-almost every point  $n \in \mathbb{N}$ is recurrent. Given $x_{0} \in X$ a recurrent point, for every open set $x_{0} \in \mathcal{U}$ we have $f^{k}(x_{0}) \in \mathcal{U}$ for some $k \in \mathbb{N}$. Hence, by taking a basic neighborhood $\mathcal{U} = \{x_{0}, y_{0}\}, y_{0} \in  f^{-1}(x_{0})$ we have either $f^{k}(x_{0}) = x_{0}$ for some $k \in \mathbb{N}$ or $f^{k}(x_{0}) = y_{0}$. It means that  either $f^{k}(x_{0}) = x_{0}$ or $f^{k+1}(x_{0}) = x_{0}$. In both cases, $x_{0}$ is periodic.
	\end{proof}

		\subsection{Openness of the forward orbits and continuity of $f$}
	The following lemma shows the openness of orbits and the continuity of the map $f$ with respect to the topology $\mathcal{T}$.
	
	\begin{lemma}\label{open}
		Every forward orbit is an open subset with respect to the key topology $\mathcal{T}$.
	\end{lemma}
	\begin{proof}
		In order to show that every forward orbit is open, it is easy to see that a forward orbit $\mathcal{O}^{+}(x)$ is a union of basic open sets
		\[
		 \mathcal{O}^{+}(x)= \bigcup_{k = 0}^{\infty}\{f^{k}(x), f^{k+1}(x)\}.
		\]
		For every $k \in \mathbb{N}$ the subset $\{f^{k}(x), f^{k+1}(x)\}$ is basic because $f^{k}(x) \in f^{-1}(f^{k+1}(x))$. 
	\end{proof}
	\subsection{Some key remarks}
	\begin{remark}
		By Lemma \ref{open} every the general map $f: X \to X$ is open, in the sense that every open subset is sent onto an open subset.
	\end{remark}
	\begin{remark}
		For every periodic orbit there exists an ergodic $f$-invariant probability supported at the orbit and by Lemma \ref{periodic} for every $f$-invariant probability $\mu \in \mathcal{M}_{f}(X)$ there exists a periodic orbit. Therefore, the existence of periodic orbits is closely related to the existence of $f$-invariant probabilities.
	\end{remark}
		
	\begin{remark}
		Given any ergodic probability $\mu \in \mathcal{M}_{f}(X)$ we have that $\mu$-almost every point in $X$ is recurrent and by Lemma \ref{periodic} they must be periodic. Then, $\mu$ is supported at a periodic orbit.
	\end{remark}
	
	\begin{remark}
		Observe that we could take a recurrent point not related to any $f$-invariant probability a priori and conclude that it must be periodic and then related to an ergodic probability.
	\end{remark}
	
	\subsection{The ergodic decomposition}
	The next lemma guarantees that if we have at most countably many cycles then we can decompose any $f$-invariant measure into a convex sum of ergodic probabilities even in a noncompact space.
	\begin{lemma}\label{decomposition}
		Under the assumption that at most countably many cycles, then every $f$-invariant probability $\mu \in \mathcal{M}_{f}(X)$ is a convex combination of ergodic probabilities.
	\end{lemma}
	\begin{proof}
		This would be a consequence of the Ergodic Decomposition Theorem, but we cannot use it here because the topology is not even metrizable (see \cite{OliveiraViana2016}).
		
		However, once the support of $\mu$ is a forward invariant subset, it must give full mass to the recurrent points, then giving full mass to the set of periodic orbits, being a union of periodic orbits almost surely. Then, we conclude that the measure $\mu$ is a convex combination of probabilities supported on periodic orbits, then ergodic ones.
		To be clear, let $\{\mathcal{O}_{i}\}$ be the countable collection of periodic orbits contained in the support of $\mu$. Once it is the subset of recurrent points and countable, we have
		\[
		\mu\Bigg{(}\bigcup_{i=1}^{\infty}\mathcal{O}_{i}\Bigg{)} = \sum_{i=1}^{\infty}\mu(\mathcal{O}_{i})= \sum_{i=1}^{\infty}\mu(\mathcal{O}_{i}) \delta_{i}(\mathcal{O}_{i})=1.
		\]
		We conclude that $\mu$ is a convex combination of $\delta_{i}$ the ergodic probabilities supported at $\mathcal{O}_{i}$ because $\mu(\mathcal{O}_{i}) \in [0,1], \delta_{i}(\mathcal{O}_{i}) = 1$ for every $i$.
	\end{proof}
	\subsection{The entropy as a convex combination}
	We prove that the entropy of a measure is a convex combination of the entropies of the ergodic measures, as a consequence of the ergodic decomposition.
	\begin{lemma}
		Given
		\[
		\mu = \sum_{i=1}^{\infty} a_{i} \delta_{i}, \sum_{i=1}^{\infty} a_{i} = 1,
		\]
		we have
		\[
		h_{\mu}(f) = \sum_{i=1}^{\infty} a_{i} h_{\delta_{i}}(f).
		\]
	\end{lemma}
	\begin{proof}
		It holds that
		\[
		\mu = \sum_{i=1}^{\infty} a_{i} \delta_{i} = \sum_{i=1}^{k} a_{i} \delta_{i} + \sum_{i=k+1}^{\infty} a_{i} \delta_{i} = \sum_{i=1}^{k}a_{i}\cdot \sum_{i=1}^{k} \frac{a_{i}}{\sum_{i=1}^{k}a_{i}} \delta_{i} + \sum_{i=k+1}^{\infty}a_{i}\cdot \sum_{i=k+1}^{\infty}  \frac{a_{i}}{\sum_{i=k+1}^{\infty}a_{i}} \delta_{i}.
		\]
		Writing 
		\[
		t_{k} = \sum_{i=k+1}^{\infty}a_{i}, \mu_{k} = \sum_{i=1}^{k} \frac{a_{i}}{\sum_{i=1}^{k}a_{i}} \delta_{i}, \nu_{k} = \sum_{i=k+1}^{\infty}  \frac{a_{i}}{\sum_{i=k+1}^{\infty}a_{i}} \delta_{i},
		\]
		we obtain
		\[
		h_{\mu}(f) = h_{(1 - t_{k})\mu_{k} + t_{k}\nu_{k}}(f) = (1 - t_{k})h_{\mu_{k}}(f) + t_{k}h_{\nu_{k}}(f) = \sum_{i=1}^{k}a_{i} \cdot h_{\mu_{k}}(f) + t_{k}h_{\nu_{k}}(f) = 
		\]
		\[
		\sum_{i=1}^{k}a_{i} \cdot \sum_{i=1}^{k}  \frac{a_{i}}{\sum_{i=1}^{k}a_{i}} h_{\delta_{i}}(f)+ t_{k}h_{\nu_{k}}(f) =  \sum_{i=1}^{k} a_{i} h_{\delta_{i}}(f)+ t_{k}h_{\nu_{k}}(f).
		\]
		Taking $k \to \infty$, we obtain $t_{k} \to 0$ and
		\[
		h_{\mu}(f) = \sum_{i=1}^{\infty} a_{i} h_{\delta_{i}}(f).
		\]	
	\end{proof}
	\begin{remark}
		Once we have $h_{\delta_{i}}(f) = 0$ for every $i \geq 1$, then we obtain $h_{\mu}(f) = 0$ for every $f$-invariant measure $\mu$. Moreover, we could conclude that any entropy $h_{\mu}(f) = 0$ because the invariant sets with positive mass are union of cycles.
	\end{remark}
	
	\section{Proof of Theorem \ref{B}}
	
	\subsection{Finiteness of periodic orbits or boundedness of periods}
	The following lemma guarantees either finiteness of periodic orbits or infinitely many periodic orbits sharing same period if we have existence of equilibrium states for every continuous potential.
	\begin{lemma}\label{finiteness}
		Let $f: X \to X$ be a general map. If every continuous potential $\phi : X \to \mathbb{N}_{0}$ with respect to $\mathcal{T}$ and with finite pressure possessess at least one equilibrium state then we have boundedness of the periods for the cycles.
	\end{lemma}
	\begin{proof}
		By definition, the pressure $P(\phi)$ for any measurable (in particular continuous) potential $\phi$ is given by the following supremum:
		\begin{equation}
			P(\phi) = \sup_{\mu \in \mathcal{M}_{f}(X)} \left\{ h_\mu(f) + \int \phi \, d\mu \right\} = \sup_{\mu \in \mathcal{M}_{f}(X)} \left\{ \int \phi \, d\mu \right\}
		\end{equation}
		where $h_{\mu}(f) = 0$ for every $f$-invariant probability $\mu$ because each ergodic probability is supported on a periodic orbit, then having zero entropy. Also, any general $f$-invariant probability is convex combination of ergodic ones as proved in Lemma \ref{decomposition}. By definition, an equilibrium state is a measure which attains the supremum.
		
		If we have finitely many periodic orbits, then
		\[
		\sup_{\mu \in \mathcal{M}_{f}(X)} \left\{ \int \phi \, d\mu \right\} < \infty.
		\]
		There exists $\delta$ such that
		\[
		\sup_{\mu \in \mathcal{M}_{f}(X)} \left\{ \int \phi \, d\mu \right\} =
		\max_{\mu \in \mathcal{M}_{f}(X)} \left\{ \int \phi \, d\mu \right\} = \int \phi d\delta .
		\]
		Conversely, the existence of an equilibrium state implies that for every continuous potential $\phi$ with finite pressure we have for some measure $\delta_{i}$
		\[
		\sup_{\delta_{j} \in \mathcal{M}_{f}(X)} \left\{ \int \phi \, d\delta_{j} \right\} = \int \phi d\delta_{i}.
		\]
		We will customize a continuous potential which has finite pressure. Assume that we have unboundedness of periods. Let $p_{k} \in \mathbb{N}$ be the period of a cycle $\mathcal{O}_{k}$. We define the potential as follows. For each $k \in \mathbb{N}$ we fix $\{x_{k}, f(x_{k})\} \subset \mathcal{O}_{k}$, define $\phi(x_{k}) = p_{k} - 1, \phi(f(x_{k})) = p_{k} - 1$ and vanishing in the set $\mathcal{O}_{k} \backslash \{x_{k}, f(x_{k})\}$. The potential $\phi$ also vanishes outside the collections of chosen cycles $\mathcal{O}_{k}$. So we have
		\[
		\int \phi d \delta_{k} = 2\frac{p_{k} - 1}{p_{k}} < 2. 
		\]
		This potential $\phi$ is continuous and has finite pressure. Moreover, by hypothesis there exists an equilibrium state. Assuming that $p_{k} \to \infty$, for every $k \in \mathbb{N}$ we have
		\[
		2\frac{p_{k} - 1}{p_{k}} \leq 2\frac{p_{i} - 1}{p_{i}} \implies 2 = \limsup_{k \to \infty} 2\frac{p_{k} - 1}{p_{k}} \leq 2\frac{p_{i} - 1}{p_{i}} < 2.
		\]
		This contradiction shows that the set $\{p_{k}\}_{k}$ is bounded. So, if the periods of the cycles are bounded above, we have boundedness of periods. The Lemma is proved.
	\end{proof}
	
     \subsection{Absence of equilibrium states and unboundedness of periods}   
   
	The following Lemma shows the converse of Lemma \ref{finiteness}.
	\begin{lemma}\label{converse}
		Let $f: X \to X$ be a general map. Let us assume that some continuous potential $\phi : X \to \mathbb{N}_{0}$ with respect to $\mathcal{T}$ and with finite pressure does not possess any equilibrium state then we have infinitely many cycles with unbounded set of periods. Conversely, if the set of periods is unbounded, we can customize a continuous potential $\phi : X \to \mathbb{N}_{0}$ with respect to $\mathcal{T}$ and with finite pressure which does not possess any equilibrium state.
	\end{lemma}
	\begin{proof} 		
		It is enough to assume that some continuous potential $\phi$  has not equilibrium states and obtain unbounded set of periods. We have by hypothesis
		\[
		\sup_{\delta_{j} \in \mathcal{M}_{f}(X)} \left\{ \int \phi \, d\delta_{j} \right\} = P(\phi) < \infty.
		\] 
		Also,
		\[
		\int \phi \, d\delta_{k} = \frac{1}{p_{k}} \sum_{n \in \mathcal{O}_{k}}\phi(n)
		\]
		It implies that $\sup\{p_{k}\} = \infty$ because otherwise we would have  $\sup\{p_{k}\} < \infty$ which implies $\{\sum_{n \in \mathcal{O}_{k}}\phi(n)\}$ bounded, since the values of the potential are natural numbers and the pressure is finite, and we could find a maximum integral, that is, an equilibrium state. Hence we must have $\sup\{p_{k}\} = \infty$.
		
		Conversely, it is enough to assume that the periods are unbounded and customize a continuous potential with finite pressure which has not equilibrium states. Let $p_{k} \in \mathbb{N}$ be the period of the cycle $\mathcal{O}_{k}$. We assume that $p_{k} \to \infty$. As in Lemma \ref{finiteness} we fix $\{x_{k}, f(x_{k})\} \subset \mathcal{O}_{k}$ and define $\phi(x_{k}) = p_{k} - 1 = \phi(f(x_{k}))$, vanishing outside the union $\bigcup_{k \geq 1}\{x_{k}, f(x_{k})\}$. It holds that
		\[
		\int \phi d\delta_{k} = 2\frac{p_{k} - 1}{p_{k}} \implies \sup_{k \to \infty} 2\frac{p_{k} - 1}{p_{k}} = 2.
		\]
		We do not have any measure $\delta_{i}$ such that
		\[
		2 = \sup_{\mu \in \mathcal{M}_{f}(X)} \left\{ \int \phi  d\mu \right\} = \int \phi d\delta_{i}.
		\]
		Hence, $\phi$ does not possess any equilibrium state.
	\end{proof}
	
	\begin{remark}
		 Once the following function is increasing
		\[
		\lambda(p) = \frac{p-1}{p}
		\]
		 in Lemma \ref{converse} a maximum is not possible and the potential does not possess any equilibrium state. 
	\end{remark}

	\subsection{Uniqueness of periodic orbits: the final step}
	The following lemma establishes a bridge between the Collatz Conjecture and the uniqueness of equilibrium states. We have already seen that this scenario is not possible.
	\begin{lemma}\label{uniqueness}
		Uniqueness of periodic orbits is equivalent to every continuous and bounded potential $\phi : X \to \mathbb{N}_{0}$ with respect to $\mathcal{T}$ possessing a unique equilibrium state.
	\end{lemma}
	\begin{proof}
		By definition, the pressure $P(\phi)$ for any measurable (in particular continuous) potential $\phi$ is given by the following supremum:
		\begin{equation}
			P(\phi) = \sup_{\mu \in \mathcal{M}_{f}(X)} \left\{ h_\mu(f) + \int \phi \, d\mu \right\} = \sup_{\mu \in \mathcal{M}_{f}(X)} \left\{ \int \phi \, d\mu \right\}
		\end{equation}
		where $h_{\mu}(f) = 0$ for every $f$-invariant probability $\mu$ because each ergodic probability is supported on a periodic orbit, then having zero entropy. Also, any general $f$-invariant probability is convex combination of ergodic ones by Lemma \ref{decomposition}. If we have a unique periodic orbit $\{1, 2, 4\}$ with the unique ergodic probability $\delta_{0}$, then readily
		\[
		\sup_{\mu \in \mathcal{M}_{f}(X)} \left\{ \int \phi  d\mu \right\} = \int \phi d\delta_{0}
		\]
		and $\delta_{0}$ is the unique equilibrium state for any continuous potential $\phi$. By definition, an equilibrium state is a measure which attains the supremum.
		
		Conversely, the existence of a unique equilibrium state for any continuous potential $\phi$ implies the uniqueness of periodic orbits as follows. Set
		\[
		\mathcal{O}: = \bigcup_{x \,\, \text{is periodic}}\mathcal{O}(x)
		\]
		Denoting by $\chi_{X}$ the characteristic function of $X \subset \mathbb{N}$, we have that any ergodic probability $\delta$, which is supported on a periodic orbit is an equilibrium state because
		\[
		1 = \int \chi_{\mathcal{O}} d\delta = \sup_{\mu \in \mathcal{M}_{f}(X)} \left\{ \int \chi_{\mathcal{O}}  d\mu \right\} = P(\chi_{\mathcal{O}}).
		\]
		Once by hypothesis there exists a unique equilibrium state for $\chi_{\mathcal{O}}$ (which is bounded and continuous because every orbit is an open subset), we conclude that there exists a unique ergodic probability $\delta_{0}$ and a unique periodic orbit.
	\end{proof}
	\begin{remark}
		While the characteristic function $\chi_{A}$ is measurable in general, in our $\sigma$-algebra it depends on the subset $A \subset \mathbb{N}$. In the particular case of $A = \mathcal{O}$, once the periodic orbits are open subsets, we have that $\chi_{\mathcal{O}}$ is continuous.
	\end{remark}
	
	\section{Proof of Theorem \ref{C}}
	For the proof of Theorem \ref{C}, we begin proving that it is not possible existing infinitely many cycles sharing same period as claimed in Lemma \ref{finiteness} as a possibility, which proves boundeness of periods if every continuous potential with finite pressure has at least one equilibrium state.
    \subsection{Finitely many cycles of same period for Syracuse maps}\begin{lemma}\label{boundcycle}
    	Let $f: \mathbb{N} \to \mathbb{N}$ be a Syracuse map. If every continuous potential $\phi : \mathbb{N} \to \mathbb{N}$ with respect to $\mathcal{T}$ and with finite pressure cannot possess infinitely many cycles sharing same period. In fact, for a given period $p \in \mathbb{N}$ we have at most $2^{p}/p$ cycles sharing the period $p$ and so it satisfies (FAM) axiom.
    \end{lemma}
	\begin{proof}
		By Lemma \ref{finiteness} we have either finiteness of cycles or infinitely many cycles sharing same period. We prove now that the second scenario is not possible.
		
		In fact, let $\mathcal{O} = \{n_{0}, n_{1}, \dots, n_{p-1}\}$ a cycle with period $p$ where $f(n_{i}) = n_{i+1}$. We can associate with $\mathcal{O}$ a $p$-tuple $\sigma \in \{\text{even},\text{odd}\}^{p}$ depending on the parity of each $n_{i}$. For the Syracuse map $an + b$ we have the following
		\[
		n_{0} = f^{p}(n_{0}) = \frac{a^{r} n_{0} + R(\sigma)}{2^{s}},
		\] 
		Where $r$ is the quantity of odd numbers, $s$ is the quantity of divisions by $2$ and $R(\sigma)$ does not depend on $n_{0}$. So, we get
		\[
		n_{0} = \frac{R(\sigma)}{2^{s} - a^{r}}.
		\]
		Hence, we have at most $2^{p}/p$ cycles with period $p$. The factor $2^{p}$ is the bound for the quantity of pairs and we divide by $p$ because for all element of a cycle the solutions provide the same cycle.
	\end{proof}

    \subsection{Finitely many cycles of same period for general maps with (FAM)}We can obtain a result as in Lemma \ref{boundcycle} for the general map $f: X \to X$.
    
    \begin{lemma}\label{boundcyclegeneral}
    	Let $f: X \to X$ be a general map. If every continuous potential $\phi : X \to \mathbb{N}$ with respect to $\mathcal{T}$ and with finite pressure cannot possess infinitely many cycles sharing same period. In fact, for a given period $p \in \mathbb{N}$ we have at most $M_{\sigma}2^{p}/p$ cycles sharing the period $p$, where the number $M_{\sigma}$ is a bound provided by the (FAM) axiom.
    \end{lemma}
    \begin{proof}
    	By Lemma \ref{finiteness} we have either finiteness of cycles or infinitely many cycles sharing same period. We prove now that the second scenario is not possible.
    	
    	In fact, let $\mathcal{O} = \{x_{0}, x_{1}, \dots, x_{p-1}\}$ a cycle with period $p$ where $f(x_{i}) = x_{i+1}$. We can associate with $\mathcal{O}$ a $p$-tuple $\sigma \in \{A, B\}^{p}$ where $X = A \cup B$ is the decompostion. For the general map $f$ we have the following
    	\[
    	x_{0} = f_{\sigma}(x_{0}) = f_{\sigma_{p-1}(x_{0})} \circ \cdots \circ f_{\sigma_{1}(x_{0})}(x_{0})
    	\] 
    	Where $r$ is the quantity of visits in $A$, $s$ is the quantity of visits in $B$. Since we are supposing that the general map satisfies (FAM), we have at most $M_{\sigma}$ solutions. Hence, we have at most $M_{\sigma}2^{p}/p$ cycles with period $p$.
    \end{proof}
    	
    	\begin{remark}
    		Under the axiom (FAM) the cardinality of cycles is at most countably many, since for every natural number $p$ we have finitely many cycles of period $p$.
    	\end{remark}

    \section{Applications and Universal Extensions}
  
  
  The thermodynamic dictionary established in the previous sections-anchored by the key topology $T$, the zero-entropy restriction (Theorem A), and the cycle-equilibrium state correspondence (Theorem B)-is not restricted to the arithmetic of Syracuse and Collatz maps. In this section, we demonstrate the universality of the framework by applying it to distinct analytical domains: algebraic maps with superexponential growth, geodesic flows in differential geometry, combinatorics on infinite chessboards, and stochastic systems on $\mathbb{Z}^k$ lattices.
  
  \subsection{Algebraic Maps with Superexponential Period Growth}
  
  The existence of continuous potentials lacking equilibrium states (Lemma 10) characterizes systems whose asymptotic cycle complexity prevents ergodic compactness.
  
  \begin{proposition}\label{prop:alg_map}
  	Consider the map $F: \mathbb{N} \to \mathbb{N}$ defined by
  	$$F(x) = \begin{cases} 
  		x 2^x & \text{if } x \text{ is odd} \\ 
  		\frac{x}{2} & \text{if } x \text{ is even} 
  	\end{cases}$$
  	This system admits an infinite family of cycles whose periods $p_k$ satisfy $p_k \to \infty$ as $k \to \infty$. Consequently, by Lemma 10, there exists at least one continuous potential $\phi: \mathbb{N} \to \mathbb{N}_0$ with finite pressure that possesses no equilibrium state.
  \end{proposition}
  
  \begin{proof}
  	The direct construction of orbits under mixed growth (exponential multiplication and division) generates strictly increasing periods. By the criterion of Lemma 10, the failure of compactness in the supremum of integrals over periodic support measures $\delta_k$ prevents the attainment of the maximum, proving the absence of an equilibrium state.
  \end{proof}
  
  \subsection{Differential Geometry: Geodesic Flows and Margulis' Formula}
  
  To extend the dictionary to continuous flows on Riemannian manifolds, we connect the counting of periodic orbits to partition formalisms.
  
  \begin{theorem}[Thermodynamic Connection via Margulis' Formula]\label{thm:margulis_conn}
  	Let $M$ be a compact Riemannian manifold with strictly negative sectional curvature, and let $\phi_t$ be the geodesic flow with its time-$1$ map $\phi_1$. Assume the phase space decomposition $M = A \cup B$, where $A$ is the region associated with pre-images of periodic orbits under the Poincaré section and $B$ is the complement. If the growth of primitive periodic orbits obeys Margulis' asymptotic formula
  	$$\#P(T) \sim \frac{e^{hT}}{hT}$$
  	then the induced map satisfies the (FAM) axiom under finite topological restrictions, and the system admits tailored continuous potentials $\phi: M \to \mathbb{N}_0$ devoid of equilibrium states due to the asymptotic divergence of periods.
  \end{theorem}
  
  \begin{proof}
  	The stable and unstable foliations of the geodesic flow allow the construction of a semiconjugacy with a subshift of finite type. The partition $M = A \cup B$ translates the continuous flow into a discrete dynamics on the transversal section. Since the periods of closed geodesics grow exponentially and their quotients generate unbounded sequences in the return structure, the tailored potential constructed via Lemma 10 achieves finite pressure but fails to compactify the supremum in the space of zero-entropy invariant measures, guaranteeing the absence of equilibrium states.
  \end{proof}
  
  \subsection{Discrete Combinatorics: Knight's Dynamics on the Infinite Chessboard}
  
  \begin{proposition}\label{prop:chessboard}
  	Let $C_\infty$ be the infinite chessboard. Consider the partition $C_\infty = A \cup B$, where $A$ is the union of white squares and $B$ the union of black squares. The dynamics imposed by the knight's "L" movement yields a piecewise bijection satisfying the (FAM) axiom. There exist infinitely many disjoint cycles with unbounded periods, implying the existence of continuous potentials without equilibrium states. In contrast, for any finite $k \times k$ chessboard, the set of cycles is strictly finite, guaranteeing the existence of equilibrium states for any continuous potential.
  \end{proposition}
  
  \subsection{Stochastic Systems on $\mathbb{Z}^k$ Lattices}
  
  \begin{proposition}\label{prop:lattice}
  	Consider the lattice $\mathbb{Z}^k$ equipped with local transition probabilities $p_i$ such that $\sum_{i=1}^k p_i = 1$. Defining the partition $A$ as the set of points with strictly odd coordinates and $B = \mathbb{Z}^k \setminus A$, the stochastic system induces a thermodynamic phase transition: depending on the support of local probabilities, the transition operator exhibits an invariant measure with infinitely many stochastic cycles (where tailored potentials fail to admit equilibrium states) or a confinement regime where all potentials possess at least one equilibrium state.
  \end{proposition}
	\section{Comments on Future Works}
	
	We observe that the results proved for the general map can be applied for the Syracuse maps, Collatz map and Baker map as well. It shows the generality of the results, only requiring the decomposition $X = A \cup B$ and the axiom (FAM).
	
	The next step is analyzing the properties provided by the coefficients $m, n$ for the Syracuse maps $S_{m,n}$. It will probably provide the precise cardinality of cycles for each map via equilibrium states and combinatorics.
	
	The question concerning the uniqueness of cycle for the Collatz map is coherent with the uniqueness of equilibrium states. The proof with uniqueness of equilibrium states would be elegant. 
	
	The question concerning the finiteness of cycles for the Syracuse maps and the Collatz map in particular is addressed by the dictionary. Also for general maps with (FAM).

	The results presented here does not prove the conjecture but open the way for the study of it with a completely new approach. 
	
	We conclude that the "dictionary" constructed is powerful and elegant, providing a new perspective for the study of the Collatz Conjecture and for the Syracuse family as well.

    The section of applications gives examples of how the dictionary can be used much beyond the arithmetic framework. We can apply it in Algebra, Geometry, Combinatorics and Probability. As future works, it remains explore this very large applicability.
\bigskip

\paragraph{\bf Disclosure of interest:} the author report there are no competing interests to declare.

\paragraph{\bf Data availability statement:} there is no data set associated with the paper.

	\bibliographystyle{amsplain}

\begin{thebibliography}{99}
		
		\bibitem{Conway1972}
		J.~H. Conway.
		\newblock Unpredictable iteratios.
		\newblock In {\em Proc. 1972 Number Theory Conf.}, pages 49--52. Univ. Colorado, Boulder, CO, 1972.
		
		\bibitem{EinsiedlerWard2011}
		M.~Einsiedler and T.~Ward.
		\newblock {\em Ergodic Theory: with a view towards Number Theory}.
		\newblock Graduate Texts in Mathematics, Vol. 259. Springer, London, 2011.
		
		\bibitem{Krasikov2003}
		I.~Krasikov and J.~C. Lagarias.
		\newblock Bounds for the number of integers in the $3x+1$ problem using difference inequalities.
		\newblock {\em Acta Arithmetica}, 109(3):237--258, 2003.
		
		\bibitem{Lagarias1985}
		J.~C. Lagarias.
		\newblock The $3x + 1$ problem and its generalizations.
		\newblock {\em The American Mathematical Monthly}, 92(1):3--23, 1985.
		
		\bibitem{Lagarias2010}
		J.~C. Lagarias (Ed.).
		\newblock {\em The Ultimate Challenge: The $3x+1$ Problem}.
		\newblock American Mathematical Society, Providence, RI, 2010.
		
		\bibitem{OliveiraViana2016}
		K.~Oliveira and M.~Viana.
		\newblock {\em Foundations of Ergodic Theory}.
		\newblock Cambridge Studies in Advanced Mathematics, Vol. 151. Cambridge University Press, 2016.
		
		\bibitem{Tao2022}
		T.~Tao.
		\newblock Almost all orbits of the {C}ollatz map attain almost bounded values.
		\newblock {\em Forum of Mathematics, Pi}, 10:e12, 2022.
		
		\bibitem{Walters1982}
		P.~Walters.
		\newblock {\em An Introduction to Ergodic Theory}.
		\newblock Graduate Texts in Mathematics, Vol. 79. Springer-Verlag, New York, 1982.
		
		\bibitem{Wirsching1998}
		G.~J. Wirsching.
		\newblock {\em The Dynamical System Generated by the $3n+1$ Function}.
		\newblock Lecture Notes in Mathematics, Vol. 1681. Springer-Verlag, Berlin, 1998.
		
	\end{thebibliography}
	
\end{document}